\documentclass[10pt, reqno]{amsart}
\usepackage[utf8x]{inputenc}
\usepackage{amssymb}
\usepackage{amsthm}
\usepackage{enumerate}

\newcommand{\al}{\alpha}
\newcommand{\ann}{\operatorname{Ann}}
\newcommand{\be}{\beta}
\newcommand{\dmssn}{\operatorname{dim}}
\newcommand{\dis}{{\operatorname{dis}}}
\newcommand{\id}{\operatorname{id}}
\newcommand{\la}{\langle}

\newcommand{\ra}{\rangle}
\newcommand{\rmu}{{r^{-1}}}
\newcommand{\smu}{{s^{-1}}}
\newcommand{\spn}{\operatorname{span}}
\newcommand{\spncl}{\overline{\spn}\, }
\newcommand{\tmu}{{t^{-1}}}

\newcommand{\vep}{\varepsilon}

\newcommand{\X}{\mathcal{X}}
\newcommand{\Y}{\mathcal{Y}}
\newcommand{\Z}{\mathcal{Z}}

\newcommand{\C}{\mathbb{C}}

\newcommand{\Kb}{\mathbb{K}}
\newcommand{\Lb}{\mathbb{L}}
\newcommand{\Mb}{\mathbb{M}}
\newcommand{\N}{\mathbb{N}}
\newcommand{\R}{\mathbb{R}}
\newcommand{\Zb}{\mathbb{Z}}

\theoremstyle{plain}
\newtheorem{theorem}{Theorem}[section]
\newtheorem{lemma}[theorem]{Lemma}
\newtheorem{corollary}[theorem]{Corollary}
\newtheorem{proposition}[theorem]{Proposition}
\theoremstyle{definition}
\newtheorem{definition}[theorem]{Definition}
\theoremstyle{remark}
\newtheorem{remark}[theorem]{Remark}
\newtheorem{example}[theorem]{Example}

\title{Construction of enveloping actions}
\author{Dami\'an Ferraro}
\date{\today}
\thanks{This work was partially supported by Mathamsud network U11-MATH05 (partially funded by ANII, Uruguay) and started during my visit to the \textit{Universidade Federal de Santa Catarina}. I thank Professors Alcides Buss and Ruy Exel for receiving me there}
\address{Departamento de Matem\'atica y Estad\'istica del Litoral, Universidad
de la Rep\'ublica, Gral. Rivera 1350. Salto. Uruguay.}
\email{dferraro@unorte.edu.uy}

\begin{document}
\begin{abstract}
  We study the problem of constructing a globalization for partial actions on *-algebras, C*-algebras and Hilbert modules.
  For the first ones we give a necessary condition for the existence of a globalization and we prove this conditions is necessary and sufficient for C*-algebras.
  Using the linking algebra of a Hilbert module we translate this condition to the realm of partial action on Hilbert modules.
\end{abstract}

\maketitle
\section{Introduction}

Among the simplest examples of partial actions on C*-algebras \cite{Ex94Circle,Ex97,Mc95} we find restrictions of actions (also called global actions) to non invariant C*-ideals \cite{Ab03}.
Many aspects of these examples, as the representations and crossed products, can be studied using techniques developed for actions on C*-algebras.
Then it is interesting to know which partial actions can be globalized, that is: described as the restriction of a global action (the globalization).
This problem was stated in \cite{Ab03}, where it was solved in case the C*-algebra is commutative.

Partial actions can be defined in other categories, as topological spaces, rings and Hilbert modules \cite{Ab03,cortes2009globalization,ExDo05,dokuchaev2007globalizations}.
Here we will work with *-algebras, C*-algebras and Hilbert modules, mainly because C*-algebras are, at once, *-algebras and Hilbert modules.

This work is organized as follows.
In the first section we study the problem of globalizing partial actions on *-algebras (*-partial actions).
Our intention is to study continuous partial action on C*-algebras (C*-partial actions) from a *-algebraic point of view, which amounts to dump all the topological structure keeping the *-algebraic properties we need.
Hence, when necessary, we make some assumptions on the *-algebras which are known to hold for C*-algebras.
Under these assumptions we give a necessary and sufficient condition for the existence of a globalization.
Then we turn to consider all the structure of C*-partial actions.
At this point we show a C*-partial action can be globalized to a C*-partial action if and only if it can be globalized to a *-partial action, independently of the topological structure.
Finally, in the last section, we give a necessary and sufficient condition for the existence of a globalization of a partial action on a Hilbert bimodule (Hb-partial action).
We specifically show that a Hb-partial action has a globalization if and only if it's linking partial action \cite{Ab03} has a globalization.

\section{Partial action on *-algebras}\label{section pa on star algebras}

Each one of the types of partial actions considered here has it's own notion of globalization.
Besides, a partial action on a C*-algebra is a partial action in a *-algebra and in a Hilbert module.
For this reason we introduce the terms ``*-partial action'', ``C*-partial action'' and ``Hb-partial action''.
In case the partial action under consideration is global we substitute the term \textit{partial} for \textit{global.}

We start by recalling some definitions and facts.

\subsection{Algebras with involution}
A *-algebra is an algebra $A$ over the complex field together with a conjugate linear function (the involution) $A\to A,\ a\mapsto a^*,$ satisfying $(ab)^*=b^*a^*$ and ${a^*}^*=a.$

Assume $A$ is a *-algebra.
By a *-ideal of $A$ we mean a subspace $I\subset A$ such that $I^*=I$ and $IA\subset A.$
A function between *-algebras, $\phi\colon A\to B,$ is a *-homomorphism if it is linear, multiplicative and $\phi(a^*)=\phi(a)^*.$
An automorphism of $A$ is a bijective *-homomorphism from $A$ to $A$ and the set of automorphism of $A$ will be denoted $Aut(A).$

An element $a\in A$ is a right (left) annihilator if $aA=\{0\}$ ($Aa=\{0\}$).
Note $a$ is a right annihilator if and only if $a^*$ is a left annihilator.
We say $A$ is non-degenerate if it does not have a right (or left) annihilator different from $0.$

\begin{remark}\label{remark equivalencies non degenerate}
  If $A$ is non-degenerate then for all $a,b\in A$ the following conditions are equivalent: (1) $a=b,$ (2) for all $c\in A,$ $ac=bc$  (3) for all $c\in A,$ $ca=cb$ and (4) for all $c,d\in A,$ $cad=cbd.$
\end{remark}

A double centraliser of $A$ is a pair $\mu=(L,R),$ where $L,R\colon A\to A$ are linear functions, $L(ab)=L(a)b,$ $R(ab)=aR(b)$ and $aL(b)=R(a)b$ (for all $a,b\in A$).
It is usual to write $\mu a :=L(a)$ and $a\mu:=R(a).$
The *-algebra structure of $M(A)$ is given by point wise vector space operations, product $(L,R)(M,S):=(L\circ M,S\circ R)$ and involution $(L,R)^*:=(L',R'),$ where $L'(a)=R(a^*)^*$ and $R'(a)=L(a^*)^*.$
The function $\tau\colon A\to M(A),$ where $\tau(a)b=ab$ and $b\tau(a)=ba,$ is *-homomorphism which is injective if and only if $A$ is non-degenerate.

\begin{remark}\label{remark multip}
  Double centralisers of *-algebras can be treated as adjointable operators of Hilbert modules.
  More precisely, if $A$ is non-degenerate then
    \begin{enumerate}
      \item Given a function $T\colon A\to A$ there exists $(L,R)\in M(A)$ such that $T=L$ if and only if there exists a function $T^*\colon A\to A,$ with $T(a)^*b=a^*T^*(b)$ for all $a,b\in A.$
      \item Given $(L_1,R_1),(L_2,R_2)\in M(A)$ we have $L_1=L_2$ if and only if $R_1=R_2.$ 
    \end{enumerate}
\end{remark}

Take any group $G$ and name $A^G$ the *-algebra of functions from $G$ to $A$ (point wise operations).
Given $t\in G$ name $\theta_t$ the automorphism of $A^G$ defined by $\theta_t(f)|_r:=f|_{rt}.$
The \textit{canonical action of $G$ on $M(A^G)$} is $\Theta\colon G\to Aut(M(A^G))$ where
$ \Theta_t(L,R)=(\theta_t\circ L\circ\theta_\tmu,\theta_t\circ R\circ \theta_\tmu).$

\subsection{Partial actions and globalizations}

Recall \cite{Mc95,Ex97} that a partial action of a group $G$ on a set $X$ is a pair $\sigma:=\left(\{X_t\}_{t\in G},\{\sigma_t\}_{t\in G}\right)$ where
\begin{enumerate}
 \item For all $t\in G,$ $X_t$ is a subset of $X$ and $\sigma_t\colon X_\tmu\to X_t$ a function.
 \item $\sigma_e$ is the identity of $X.$
 \item For all $s,t\in G,$ if $x\in X_\tmu$ and $\sigma_t(x)\in X_\smu,$ then $x\in X_{\tmu\smu}$ and $\sigma_s(\sigma_t(x))=\sigma_{st}(x).$
\end{enumerate}
In case $X_t=X$ for all $t\in G,$ $\sigma$ is said to be global and it is just a common action of $G$ on $X.$

Let us assume that $\sigma$ and $\tau$ are partial actions of $G$ on the sets $X$ and $Y,$ respectively.
We say $f\colon \sigma\to \tau$ is a \textit{morphism of partial actions on sets} if $f$ is a function from $X$ to $Y$ and for all $t\in G:$ $f(X_t)\subset Y_t$ and $f(\sigma_t(x))=\tau_t(f(x))$ (for all $x\in X_\tmu$).
The identity morphism associated to $\sigma$ is $id_\sigma:=id_X$ and the composition of morphism is just the composition of functions.

With $\sigma$ as before take a set $U\subset X.$
Given $t\in G$ define $U_t:=U\cap \sigma_t(X_\tmu\cap U).$
It is obvious that $U_\tmu \subset X_\tmu$ and $\sigma_t(U_\tmu)\subset U_t;$ hence it makes sense to define $\kappa_t\colon U_\tmu\to U_t$ as $\kappa_t(x)=\sigma_t(x).$
Straightforward arguments imply the \textit{restriction of $\sigma$ to $U$}, defined as $\sigma|_U:=\left(\{U_t\}_{t\in G},\{\kappa_t\}_{t\in G}\right),$ is partial action of $G$ on $U.$
From \cite{Ab03} we know every partial action on a set is isomorphic to the restriction of a global action on a set.

A set $U\subset X$ is said to be $\sigma-$invariant if for all $t\in G,$ $\sigma_t(X_\tmu\cap U)\subset U.$

\begin{remark}\label{remark double restriction}
  If $U\subset V\subset X$ then $\sigma|_U=\sigma|_V|_U.$
  Besides, if $\sigma$ is global, $\sigma|_U$ is global if and only if $U$ is $\sigma-$invariant.

\end{remark}

\begin{definition}
  A *-partial action of the (discrete) group $G$ on the *-algebra $A$ is a set theoretic partial action of $G$ on $A,$ $\al =\left(\{A_t\}_{t\in G},\{\al_t\}_{t\in G}\right),$ such that $A_t$ is a *-ideal of $A$ and $\al_t$ a *-homomorphism (for all $t\in G$).
\end{definition}

\begin{example}\label{example restriction of pa}
  Let $\al$ be a *-global action of $G$ on the *-algebra $A$ and let $I$ be a *-ideal of $A.$
  The restriction of $\al$ to $I,$ $\al|_I:=\left(\{I\cap \al_t(I)\}_{t\in G},\{\al_t|_{I\cap \al_t(I)}\}_{t\in G}\right),$ is a *-partial action of $G$ on $I.$
\end{example}

\begin{example}\label{example pa given by an ideal}
  Let $A$ be a $*-$algebra and $I$ a *-ideal of it.
  The *-partial action of $\Zb_2=\{0,1\}$ on $A$ determined by $I$ is $\alpha^{AI}:=\left(\{A_0,A_1\},\{\al_0,\al_1\}\right)$ where $A_0=A,$ $A_1=I,$ $\al_0=id_A$ and $\al_1=id_I.$
\end{example}

Morphisms of *-partial actions are just morphisms of partial actions on set which are also *-homomorphisms.
The composition and identity are the natural ones.

\begin{definition}[\cite{Ab03,ExDo05}]\label{definition star globalization}
  A *-globalization the *-partial action $\al$ (of $G$ on $A$) is a 4-tuple $\Xi=(B,\beta,I,\iota)$ where:
  $B$ is a *-algebra, $\be$ is a *-global action of $G$ on $B,$ $I$ is a *-ideal of $B$ and $\iota\colon \al \to \be|_I$ is an isomorphism.
  
  We say $\Xi$ is minimal if $[I]:=\spn \{ \be_t(I)\colon t\in G \}$ equals $B$ and, for convenience, we say $\Xi$ is non degenerate if $B$ is non degenerate.
\end{definition}

Do not confuse our concept of minimality with topological minimality, here $\be$ may have many open invariant sets\footnote{Note that there is no topology involved.}.
Note that $\Xi$ is minimal iff the only $\be-$invariant *-ideal of $B$ containing $I$ is $B$ itself.
In case $\Xi$ is not minimal, $([I],\beta|_{[I]},I,\iota)$ is minimal because $[I]$ is invariant and $\be$ global (Remark \ref{remark double restriction}).

\begin{example}\label{example extremely non degenerate algebras}
  Given a complex vector space $V$ and a conjugate linear bijection $T\colon V\to V$ with $T^2=id_V,$ let $V^T$ be the *-algebra obtained by considering on $V$ the null product ($uv=0$) and $T$ as involution.
  Here we consider the entry wise conjugation $T\colon \C^n\to \C^n.$
  Let $\al$ be the partial action of $\Zb_3=\{0,1,2\}$ (with additive notation) on $\C^T$ such that $\al_0 = id_{\C^T}$ and $\al_1=\al_2=id_{\{0\}}.$
  There are two completely different globalizations for $\al,$ $({\C^2}^T,\beta,I,\iota)$ and $({\C^3}^T,\gamma,J,\kappa),$ that we now describe.
  Just set $I=\{(u,0)\colon u\in \C\},$ $J=\{(u,0,0)\colon u\in \C\},$ $\iota(u)=(u,0),$ $\kappa(u)=(u,0,0),$ $$\be_1(u,v)=-\frac{1}{2}(u+\sqrt{3}v,-\sqrt{3}u+v)\mbox{ and } \gamma_1(u,v,w)=(w,u,v).$$
  Note $\be_1$ is the rotation by angle $2\pi/3$ and a simple dimension argument implies that $\be$ is not isomorphic to $\gamma.$
  In fact it can be shown that these are the unique globalizations of $\al$ (up to isomorphism).
\end{example}

\begin{example}\label{example globalize al2ai}
  Consider the partial action of Example \ref{example pa given by an ideal}.
  Assume there exists a *-ideal $J$ of $A$ such that $A=I\oplus J.$
  Now form the (external) direct sum $B:=I\oplus J \oplus J,$ considered as a *-algebra with entry wise operations.
  Let $\beta$ be the action of $\Zb_2$ on $B$ given by $\be_1(a\oplus b\oplus c)=a\oplus c\oplus b$ and $\iota\colon A\to B$ the unique linear map such that $\iota(a+b)=a\oplus b,$ for all $a\in I$ and $b\in J.$
  Then $(B,\beta,I\oplus J\oplus 0,\iota)$ is a minimal *-globalization of $\al^{AI}.$
\end{example}

The existence of the direct complement $J$ is not necessary for the existence of a globalization, as we show with the next Example.
Moreover, with it we also show that there are *-partial actions on non degenerate *-algebras with a *-globalization but without a non degenerate *-globalization.

\begin{example}\label{example pa on matrices}
  Let $U\in \Mb_4(\C)$ be the matrix corresponding to the permutation $(1\ 4)(2\ 3)$ (written as a product of cycles).
  Define an involution in $\Mb_4(\C)$ by the formula $a^*:=u \overline{a}^tu,$ where $a\mapsto \overline{a}$ is the entry wise complex conjugation and $a\mapsto a^t$ is the usual matrix transposition\footnote{With this structure $\Mb_4(\C)$ is not a C*-algebra because $a^*a=0$ if $a_{i,j}=\delta_1(i)\delta_1(j).$}.
  Note that $a^*$ is obtained from $\overline{a}$ by performing a symmetry with respect to the anti-diagonal\footnote{From lower left corner to upper right corner.}.

  Name $A$ the $*-$subalgebra of $\Mb_4(\C)$ formed by the matrices of the form
  $$ \left(\begin{array}{cccc}
          a_{11} & a_{12} & a_{13} & a_{14}\\
              0   &    0   & a_{23} & a_{24}\\
              0   &    0   &   0    & a_{34}\\
              0   &    0   &   0    & a_{44}
          \end{array}\right) $$
  With $I:=\{a\in A\colon a_{23}=0\}$ we have $I=I^*$ and $AA\subset I.$
  Then $I$ is a *-ideal of $A$ and it can be shown that $A$ and $I$ are non-degenerate.
  Now let $\al^{AI}$ be the *-partial action described in Example \ref{example pa given by an ideal}.

  Assume $\Xi=(B,\be,J,\iota)$ is a *-globalization of $\al^{AJ}.$
  For convenience we think $J=A,$ $\iota=\id_A$ and $\al=\be|_J.$
  Since $B=A+\be(A)$ and $\dmssn (A\cap \be(A))=\dmssn(I)=7,$ we know $\dmssn(B)=\dmssn(A)+1=9.$
  Let $u\in A$ be the matrix with $1$ in the entry $2-3$ and $0$ elsewhere.
  Then $A=I\oplus \C u$ and $\be_1(u)\notin A$ because $u\notin I=A\cap \be_1(A).$
  Moreover, $B=A\oplus \C\be_1(u).$
  As a vector space $B$ is isomorphic to the external direct sum $A\oplus \C $ and, with this notation, the involution of $B$ is $(a\oplus \lambda )^* = a^* \oplus  \overline{\lambda}$ (because $u^* = u$).

  To describe the product of $B$ we compute
  $$(a+\lambda \be_1(u))(b+\mu \be_1(u)) = ab + \mu a\be_1(u) + \lambda\be_1(u)b + \lambda\mu\be_1(u^2).$$
  Firstly note that $u^2=0,$ so $\be_1(u^2)=0.$
  Secondly, $a\be_1(u)\in I$ and for all $c\in I$ we have $(a\be_1(u)-au)c=a\be_1(u)\al_1(c)-auc=a\al_1(uc)-auc=0.$
  Since $I$ is non-degenerate and $au\in I,$ $a\be_1(u)=au.$
  In the same way we deduce that $\be_1(u)b=ub.$
  Now the product of $B$ takes the form
  $$ (a+\lambda \be_1(u))(b+\mu \be_1(u)) = ab + \mu au + \lambda ub. $$
  Then the formula for the product of $A\oplus \C$ should be 
  $$(a\oplus \lambda)(b\oplus \mu):= ab + \mu au + \lambda ub \oplus 0.$$
  In fact $A\oplus \C$ is a *-algebra with this product and the involution described in the previous paragraph.
  Moreover, it is degenerate because $u\oplus -1$ is a bilateral annihilator.

  The final step is to determine $\be,$ which amounts to give an expression for $\be_1.$
  Assume $b=a+\lambda u + \mu \be_1(u),$ with $a\in I$ and $\lambda,\mu\in \C.$
  Then $\be_1(b)=\be_1(a)+\lambda \be_1(u)+\mu u=a+\mu u+\lambda \be_1(u).$
  In terms of $A\oplus \C$ $\be_1$ should be given by $\be_1((a+\lambda u)\oplus \mu ) = (a\oplus \mu u)\oplus \lambda.$
  Now the reader can verify that $\be_1$ is in fact a *-homomorphism of $A\oplus \C$ with $\be_1^2=\id_{A\oplus \C}$ and $\be_1|_A = \al.$
  If we set $\iota\colon A\to A\oplus \C$ as the natural inclusion then $(A\oplus \C,\be,A\oplus 0,\iota)$ is a minimal *-globalization of $\al.$
  Furthermore, it is essentially the unique *-globalization of $\al.$
\end{example}

The algebras from Example \ref{example extremely non degenerate algebras} are those in which all the elements are annihilators.
Our source of inspiration are partial action on C*-algebras and these algebras are non-degenerate because on such algebras the identity $aa^*=0$ implies $a=0.$
Then we will assume our *-algebras are non degenerate when needed.

Now our idea is to take a *-partial action with a *-globalization and try to construct another *-globalization using just the *-partial action.
In this way we will be sure that the *-partial action has a *-globalization every time the construction can be performed.
This new *-globalization may be completely different from the initial one.

Fix, for the rest of this section, a $*-$partial action $\al$ of $G$ on $A$ and a minimal *-globalization $\Xi=(B,\beta,I,\iota)$ of $\al.$
The \textit{canonical morphism} associated to $\Xi,$ denoted $\pi$ or $\pi_\Xi$ in case is necessary to mention $\Xi,$ is the *-homomorphism $\pi\colon B\to M(A^G),$ $\pi(b)f|_r=\iota^{-1}(\be_r(b)\iota(f|_r)).$
A simple computation shows that $\pi \colon \be\to \Theta$ is a morphism, where $\Theta$ is the canonical action of $G$ on $M(A^G).$
Then $\pi(B)$ is a $\Theta$-invariant *-subalgebra and $\pi(A)$ is a *-ideal of it.
Moreover, $\Theta|_{\pi(B)}$ is a *-global action and $\Theta|_{\pi(B)}|_{\pi(I)}=\Theta|_{\pi(I)}.$
The quadruple $\Xi_\Theta:=(\pi(B),\Theta|_{\pi(B)},\pi(I),\pi\circ \iota)$ is a minimal *-globalization of $\alpha$ if and only if $\pi\circ\iota$ is an injective function and $\pi(I_t) = \pi(I)\cap \Theta_t(\pi(I))$ (for all $t\in G$).
Unfortunately these two conditions seem to be unrelated in general.

\begin{example}\label{example canonical morphism non injective}
  Let $\al^{AI}$ be the *-partial action from Example \ref{example pa on matrices} and $(A\oplus\C,\be,A\oplus 0,\iota)$ it's *-globalization.
  Then $A^{\Zb_2} = A\oplus A$ and $\pi(a\oplus \lambda)(b\oplus c)=(a+\lambda u)b\oplus (a+\lambda u)c.$
  Since $\pi(a\oplus \lambda) = \pi(a+\lambda u\oplus 0),$ we have $\pi(A\oplus \C)=\pi(A)$ and $\pi$ is not injective.
  But $\pi\circ \iota$ is injective because $A$ is non degenerate.
  In this case the kernel of $\pi,$ $\ker(\pi),$ is exactly the space generated by $u\oplus -1,$ which is the set of right annihilators of $A\oplus \C,$ $\ann_R(A\oplus \C).$
  This last condition is not accidental, as we now show.
\end{example}

\begin{lemma}\label{lemma injective non degenerate}
 Let $\Xi=(B,\be,I,\iota)$ be a *-globalization of the *-partial action $\al$ (of $G$ on $A$).
 If $\pi$ is the canonical morphism associated to $\Xi$ then $\ker(\pi)=\ann_R(B).$
 Moreover, if $A$ is non-degenerate then $\pi\circ\iota$ is injective and $\pi(B)$ is non degenerate.
\end{lemma}
\begin{proof}
  Assume $b\in \ker(\pi).$
  Given $a\in A$ and $r\in G$ denote $b\delta_r$ the element of $A^G$ taking the value $b$ at $r$ and $0$ elsewhere.
  Then $0=\iota(\pi(b)a\delta_r|_r)=\be_r(b)\iota(a).$
  This implies $bB=\spn b\be_r(\iota(A))=0,$ so $b\in \ann_R(B).$
  For the converse assume $b\in \ann_R(B)$ and take $f\in A^G.$
  Then for all $r\in G$ we have $\be_r(b)\in \ann_R(B),$ thus $\pi(b)f|_r=\iota^{-1}( \be_r(b)\iota(f|_r) )=0$ and $b\in \ker(\pi).$
  
  In case $A$ is non-degenerate and $a\in \ker(\pi\circ \iota),$ we have $\iota(a)\in \ann_R(B).$
  This implies $a\in \ann_R(A)$ and, so, $a=0.$
  
  Finally assume $\pi(b)\in \ann_R(\pi(B)).$
  Since $\pi(B)$ is $\Theta-$invariant, for all $r\in G$ we have $\pi(\be_r(b))=\Theta_r(\pi(b))\in \ann_R(\pi(B)).$
  Then, for all $r\in G$ and $c,d\in A,$ $0=\pi(\be_r(b))\pi(\iota(c))d\delta_e|_e =\iota^{-1}(\be_r(b)  \iota(c)\iota(d) ). $
  This implies $0=\be_r(b)\iota(c)$ for all $r\in G$ and $c\in A.$ Using the definition of $\pi(b)$ we conclude that $\pi(b)=0.$
\end{proof}

To give a sufficient condition for $\pi\circ\iota$ to be an isomorphism we introduce assimilative ideals.
A subset $S$ of the *-algebra $C$ is \textit{assimilative} (in $C$) if given $c\in C$ such that $cC\subset S,$ then $c\in J.$
It is immediate that every subset of a unital *-algebra is assimilative.
Furthermore, a *-ideal $J$ of $C$ is assimilative iff $C/J$ is non-degenerate.

\begin{remark}
  Recall \cite{dokuchaev2007globalizations} that a ring $R$ is $s$-unital if for all $a\in R$ there exists $b\in R$ such that $ba=a.$
  Evidently every $s$-unital *-algebra is non-degenerate and any subset of a $s$-unital *-algebra is assimilative in the algebra.
  Any sum of $s$-unital ideals is $s$-unital \cite[Remark 2.5]{dokuchaev2007globalizations}.  
\end{remark}

The next Proposition resumes all we have to say about uniqueness and non degeneracy of globalizations.

\begin{proposition}\label{proposition of uniqueness}
  Let $\al$ be a *-partial action of $G$ on $A$ with $A$ non-degenerate and $A_t$ assimilative in $A,$ for all $t\in G.$
  Then
  \begin{enumerate}[(1)]
    \item $\al$ has a *-globalization if and only if $\al$ has a non-degenerate and minimal *-globalization.
    \item Given minimal *-globalizations of $\al,$ $\Xi=(B,\beta,I,\iota)$ and $\Sigma =(C,\gamma,J,\kappa),$ with $\Sigma$ non-degenerate, there exists a unique morphism ${}_{\Xi}\rho_\Sigma\colon \beta\to \gamma$ such that ${}_{\Xi}\rho_\Sigma\circ \iota = \kappa.$
    Moreover, ${}_{\Xi}\rho_\Sigma$ is an isomorphism iff $\Xi$ is non-degenerate iff $\pi_\Xi$ is injective.
    \item In case for all $n\in \N$ and $t_1,\ldots,t_n\in G$ the ideal $A_{t_1}+\cdots+A_{t_n}$ is assimilative in $A,$ every minimal *-globalization of $\al$ is non-degenerate.
  \end{enumerate} 
\end{proposition}
\begin{proof}
  Assume $\al$ has a $*-$globalization $\Xi=(B,\beta,I,\iota),$ which we can assume is minimal (see the comments after Definition \ref{definition star globalization}).
  Let $\pi$ be the canonical morphism associated to $\Xi.$
  We claim that $\pi(I_t)=\pi(I)\cap \Theta_t(\pi(I)).$ 
  Indeed, note that $\pi(I_t)=\pi(I\cap \be_t(I))\subset \pi(I)\cap \Theta_t(\pi(I)).$
  For the converse inclusion take $a,b\in I$ and $t\in G$ such that $\pi(a)=\Theta_t(\pi(b)).$
  Given $c\in A,$ $a\iota(c)=\iota(\pi(a)\delta_e^c|_e) =\iota(\Theta_t(\pi(a))\delta_e^c |_e) =\be_t(a)\iota(c)\in I_t.$
  As $I_t$ is assimilative in $I,$ $a\in I_t$ and $\pi(I_t)=\pi(I)\cap \Theta_t(\pi(I)).$
  
  It is straightforward to show that $\pi\circ \iota(A)=\pi(I)$ is *-ideal of $\pi(B),$ $[\pi(I)]=\pi(B).$
  From the last paragraph and Remark \ref{remark equivalencies non degenerate} we know $\pi(B)$ is non-degenerate and $\pi\circ \iota\colon \al\to \Theta|_{\pi(I)}$ is an isomorphism.
  As $\Theta|_{\pi(I)}=\Theta|_{\pi(B)}|_{\pi(I)},$ $(\pi(B),\Theta|_{\pi(B)},\pi(I),\pi\circ \iota)$ is a non-degenerate and minimal *-globalization of $\al.$
  Thus we have proved (1).
  
  Lets prove the existence claim of statement (2).
  Let $\pi_\Xi$ be the canonical morphism associated to $\Xi.$
  To show that the image of $\pi_\Xi,$ $Im(\pi_\Xi),$ is $Im(\pi_\Sigma)$ note that given $(t,a,b)\in G,$ the element $u^\al_t(a,b):=\iota^{-1}(\beta_t(\iota(a))\iota(b))$ is the unique element of $A_t$ such that $\al_t(c)u^\al_t(a,b) = \al_t(ca)b$ (for all $c\in A_\tmu$).
  Indeed, as $\iota\colon \al\to \be|_I$ is an isomorphism and $\iota(A_t)=I\cap\be_t(I)$ an ideal of $B:$
  \begin{align*}
   \al_t(c)u^\al_t(a,b) 
    & = \iota^{-1}( \be_t(\iota(c))\be_t(\iota(a))\iota(b) ) = \iota^{-1}( \be_t(\iota(c))\be_t(\iota(a))\iota(b) )\\
    &=  \iota^{-1}( \be_t(\iota(ca))\iota(b) ) = \al_t(ca)b.
  \end{align*}
  Uniqueness follows from the non degeneracy of $A_t=\al_t(A_\tmu).$
  Then uniqueness implies $\iota^{-1}(\beta_t(\iota(a))\iota(b)) = u^\al_t(a,b) =\kappa^{-1}(\gamma_t(\kappa(a))\kappa(b)).$
  This last equality implies that $\pi_\Xi\circ \iota = \pi_\Sigma\circ \kappa.$
  Then $\pi_\Xi(I) = \pi_\Xi\circ \iota(A)  =\pi_\Sigma\circ \kappa(A) = \pi_\Sigma(J)$ and $Im(\pi_\Xi) = [\pi_\Xi(I)] = [\pi_\Sigma(J)]=Im(\pi_\Sigma).$
  
  From Lemma \ref{lemma injective non degenerate} we know $\ker(\pi_\Sigma)=\ann_R(C).$
  Since $C$ is non-degenerate this implies $\pi_\Sigma$ is injective.
  Now define ${}_\Sigma\rho_\Xi:=\pi_\Sigma^{-1}\circ \pi_\Xi.$
  In case $\Xi$ is non degenerate, $\pi_\Xi$ is injective and $\pi_\Xi^{-1}\circ \pi_\Sigma$ is the inverse of ${}_\Sigma\rho_\Xi.$
  To close the circle note that in case ${}_\Sigma\rho_\Xi$ is an isomorphism, $\Xi$ is non degenerate because $\Sigma$ is non degenerate.  
  Uniqueness of ${}_\Sigma\rho_\Xi$ follows from tree facts: it is a morphism of *-global actions, ${}_\Sigma\rho_\Xi|_I =\kappa\circ\iota^{-1}$ and $\spn \{\be_t(I)\colon t\in G\}=B.$
  
  To prove (3) assume $\Xi$ is minimal, we will show that $\Xi$ is non-degenerate.
  From the hypothesis we get that, for all $n\in \N$ and $t_1,\ldots,t_n\in G,$ the ideal $I_{t_1}+\cdots +I_{t_n}$ is assimilative in $I.$
  Assume $b\in B$ satisfies $Bb=\{0\}.$
  Given that $[I]=B,$ there exits $n\in \N,$ $t_1,\ldots,t_n\in G$ and $b_1,\ldots,b_n\in I$ such that $b=\sum_{j=1}^n \be_{t_j}(b_j).$
  We show $b=0$ by induction in $n.$
  If $n=1$ it follows that $Ib_1=\{0\}$ and this implies $b=0$ because $I$ is non-degenerate.
  For $n>1$ we have $\be_r(b)c = \be_r(b\be_\rmu(c)) = 0,$ for all $c\in I$ and $r\in G.$
  Then $\sum_{j=1}^n \be_{r t_j}(b_j)c=0.$
  With $r=t_{n}^{-1}$ we get $b_n c = -\sum_{j=1}^{n-1} \be_{t_n^{-1} t_j}(b_j)c \in \sum_{j=1}^{n-1} I_{t_n^{-1} t_j},$ for all $c\in I.$
  Thus there exists $b'_j\in  I_{t_n^{-1} t_j}$ ($j=1,\ldots,n-1$) such that $b_n = \sum_{j=1}^{n-1} b'_j.$
  Besides, $\be_{t_n}(b'_j) \in \be_{t_n}(I_{t_n^{-1} t_j})\subset \be_{t_n}\circ \be_{t_n^{-1}t_j}(I)=\be_{t_j}(I),$ so there exists $b''_j\in I$ such that $\be_{t_n}(b_n) = \sum_{j=1}^{n-1} \be_{t_j}(b''_j).$
  This implies $b= \sum_{j=1}^{n-1} \be_{t_j}(b_j + b''_j).$
  By induction it follows that $b=0.$
\end{proof}

Before giving a sufficient condition for the existence of a *-globalization we prove the following

\begin{lemma}\label{lemma equality}
  Let $A$ be a *-algebra, $I$ and $J$ *-ideals of $A,$ $a\in I$ and $b\in J$ such that: $I\cap J$ and $I+J$ are non-degenerate, $Ia\subset J,$ $Jb\subset I$ and $ca=cb$ for all $c\in I\cap J.$
  Then $a=b.$
\end{lemma}
\begin{proof}
For all $z\in I\cap J$ and $x\in I\cup J$ we have $xa,xb\in I\cap J$ and $zxa = zxb.$
Thus, for all $x\in I+J,$ $xa=xb$ and this implies $a=b.$
\end{proof}

The next result is a combination of \cite[Theorem 4.5]{ExDo05} and \cite[Theorem 3.1]{dokuchaev2007globalizations}.

\begin{theorem}\label{theorem enveloping star algebra}
    Let $\al=\left(\{A_t\}_{t\in G},\{\al_t\}_{t\in G}\right)$ be a *-partial action of $G$ on $A$ and consider the conditions
    \begin{enumerate}[(1)]
     \item $\al$ has a *-globalization.
     \item For all $(t,a,b)\in G\times A\times A$ there exists $u\in A_t$ such that, for all $c\in A_\tmu,$ $\al_t(c)u=\al_t(ca)b.$
    \end{enumerate}
    Then (1) implies (2).
    In case for all $s,t\in G$ the *-subalgebras $A_t+A_s$ and $A_t\cap A_s$ are non-degenerate and $A_t$ is assimilative in $A,$ (2) implies (1).
    Moreover, in this last case $\al$ has a non-degenerate and minimal *-globalization.
\end{theorem}
\begin{proof}
    The implication (1) $\Rightarrow$ (2) is part of the proof of Proposition \ref{proposition of uniqueness}.
    For the converse note that the element $u$ in (2) is uniquely determined by $(t,a,b)$ because $A_t$ is non-degenerate and $\al_t(A_\tmu)=A_t.$
    The expression $u_t(a,b)$ will be used to denote the element $u$ given for $(t,a,b)$ in (2). 
    
    Let $A^G$ be the *-algebra of function from $G$ to $A$ with point wise operations and $M(A^G)$ it's multiplier algebra.
    Note $A=A_e$ is non-degenerate, so $A^G$ is non-degenerate and we can appeal to Remark \ref{remark multip}, what we do without explicit mention.
    We claim that there exists a unique injective homomorphism of *-algebras $\rho\colon A\to M(A^G)$ such that $\rho(a)f|_t = u_t(a,f|_t).$
    
    Given $a\in A,$ to show $\rho(a)$ is a multiplier with adjoint $\rho(a^*)$ it suffices to show that $[\rho(a)f]^*g=f^*[\rho(a^*)g],$ for all $f,g\in A^G.$
    The equality holds because for all $f,g\in A^G, $ $t\in G$ and $c,d\in A_t$ we have $[\rho(a)f]^*g|_t,f^* [\rho(a^*)g]|_t\in A_t$ and
    \begin{align*}
	d [\rho(a)f]^*g|_t c
	    & = d u_t(a,f|_t)^*g|_tc
		= d [ \al_t(  \al_\tmu( c^*g^*|_t) a ) f|_t  ]^* \\
	    & = d f^*|_t \al_t( a^* \al_\tmu( g|_t c ) )
		= \al_t( \al_\tmu(d f^*|_t)  a^* \al_\tmu( g|_t c ) )\\
	    & = \al_t( \al_\tmu(d f^*|_t)  a^*)  g|_t c
		= d f^* [\rho(a^*)g]|_t c.
    \end{align*}
    
    The uniqueness of elements $u_t(a,b)$ implies $a\mapsto u_t(a,b)$ is linear, so $\rho$ is linear.
    To show $\rho$ is multiplicative it suffices to prove, for all $a,b,c\in A$ and $t\in G,$ the identity $u_t(a,u_t(b,c))=u_t(ab,c).$
    Note $u_t(a,u_t(b,c)),u_t(ab,c)\in A_t$ and for all $d\in A_t:$
    $$
	d u_t(a,u_t(b,c))
	     = \al_t(\al_\tmu(d)a)u_t(b,c)
		= \al_t(\al_\tmu(d)ab)c
		= du_t(ab,c).
    $$
    Thus $u_t(a,u_t(b,c))=u_t(ab,c).$
    
    Assume $\rho(a)=0.$
    Given $b\in A$ let $\delta_e^b\in A^G$ be the function taking the value $b$ on $e$ and $0$ elsewhere.
    Then $0=\rho(a)\delta_e^b|_e = u_e(a,b).$
    Given $c\in A=A_e$ we have $0=cu_e(a,b) =\al_e(ca)b=cab.$
    As $A$ is non-degenerate, $a=0$ and $\rho$ is injective.
    
    With $\Theta$ being the canonical action of $G$ on $M(A^G)$ set     $$B:=\spn \{\Theta_t(\rho(A))\colon t\in G\}.$$
    So $B$ is a $\Theta$-invariant *-sub algebra of $M(A^G)$ and the restriction of $\Theta$ to $B,$ $\beta,$ is a *-global action.
    The proof will be completed once we show $\rho(A)$ is an ideal of $B,$ that $\rho$ is an isomorphism between $\al$ and $\be|_{\rho(A)}$ and that $B$ is non-degenerate.
    
    Lets show that $\be_t(\rho(A))\rho(A)\subset \rho(A_t),$ for all $t\in G.$
    It suffices to show that $\be_t(\rho(a))\rho(b)=\rho(u_t(a,b)).$
    From the definitions of $\be$ and $\rho$ we obtain 
    $$\be_t(\rho(a))\rho(b)f|_r = u_{rt}(a,u_r(b,f|_r))\mbox{ and } \rho(u_t(a,b))f|_r = u_r(u_t(a,b),f|_r).$$
    Putting $c:=f|_r$ it suffices to show $ u_{rt}(a,u_r(b,c))= u_r(u_t(a,b),c).$
    
    Note $u_{rt}(a,u_r(b,c))\in A_{rt},$ $u_r(u_t(a,b),c)\in A_r,$
    $$ A_{rt} u_{rt}(a,u_r(b,c) = \al_{rt}\left( A_{\tmu\rmu} a \right)u_r(b,c) \in A_{rt}\cap A_r \mbox{ and} $$
    $$ A_ru_r(u_t(a,b),c) = \al_r\left( A_\rmu u_t(a,b) \right)c\in \al_r(A_\rmu\cap A_t)=A_r\cap A_{rt}. $$
    Besides, for $z\in A_{rt}\cap A_r:$
    \begin{align*}
      zu_{rt}(a,u_r(b,c))
	& = \al_{rt}(\al_{\tmu\rmu}(z)a)u_r(b,c)
	   = \al_r\left( \al_t( \al_{\tmu\rmu}(z)a )b \right)c\\
	& = \al_r\left( \al_\rmu(z) u_t(a,b) \right)c
	  = zu_r\left( u_t(a,b),c \right).
    \end{align*}
    Then Lemma \ref{lemma equality} implies $u_{rt}(a,u_r(b,c))= u_r(u_t(a,b),c).$
    
    Now the inclusion $\be_t(\rho(A))\rho(A)\subset \rho(A_t)$ (valid $\forall \ t\in G$) implies $\rho(A)$ is an ideal of $B$ because
    $$B\rho(A)=\sum_{t\in G} \be_t(\rho(A))\rho(A)\subset \rho(A)\mbox{ and } \rho(A)B = (B\rho(A))^*\subset \rho(A).$$
    
    Take $t\in G$ and $a\in A_\tmu.$ To prove $\rho(\al_t(a))=\be_t(\rho(a))$ is equivalent to prove $u_r(\al_t(a),b)=u_{rt}(a,b),$ for all $b\in A$ and $r\in G.$
    As before note that $u_r(\al_t(a),b)\in A_r,$ $u_{rt}(a,b)\in A_{rt},$
    $$ A_r u_r(\al_t(a),b)=\al_r\left( A_\rmu \al_t(a)\right)b \in \al_r(A_\rmu\cap A_t)=A_r\cap A_{rt}\mbox{ and} $$
    $$ A_{rt}u_{rt}(a,b)=\al_{rt}\left(A_{\tmu\rmu} a  \right)b\in \al_{rt}\left(A_{\tmu\rmu} A_\tmu  \right) = A_{rt}\cap A_r. $$
    For $z\in A_r\cap A_{rt}:$
    \begin{align*}
      zu_r(\al_t(a),b)
	& = \al_r(\al_\rmu(z)\al_t(a) )b
	  = \al_r(\al_t(\al_{\tmu\rmu}(z)a) )b\\
	& = \al_{rt}(\al_{\tmu\rmu}(z)a)b
	  = zu_{rt}(a,b).
    \end{align*}
    Then $u_r(\al_t(a),b)=u_{rt}(a,b)$ by Lemma \ref{lemma equality}.
    
    To finish the proof all we need to show is that $\rho(A_t)=\rho(A)\cap \be_t(\rho(A)).$
    Note $\rho(A_t)=\rho(\al_t(A_\tmu)) =\be_t(\rho(A_\tmu))\subset \rho(A)\cap \be_t(\rho(A)).$
    Now take $T\in \rho(A)\cap \be_t(\rho(A)).$
    Then there exist $a,b\in A$ such that $T=\rho(a)=\be_t(\rho(b)).$
    Thus, for all $f\in A^G$ and $r\in G:$
    $$ u_r(a,f|_r) =\rho(a)f|_r = \be_t(\rho(a))f|_r = u_{rt}(a,f|_r).$$
    By replacing $f$ with $\delta_e^b$ and $r$ with $e$ we obtain $u_e(a,b)=u_t(a,b),$ for all $b\in A.$
    But $u_e(a,b)=ab,$ then $aA\subset A_t$ and this implies $a\in A_t.$
    
    Up to this point we have shown that $\Xi:=(B,\Theta|_B,\rho(A),\rho)$ is a minimal *-globalization of $\al.$
    Note that the canonical morphism associated to $\Xi,$ $\pi_\Xi\colon B\to M(A^G),$ is just the natural inclusion $B\subset M(A^G)$ because $\pi_\Xi\circ \rho =\rho.$
    Then Lemma \ref{lemma injective non degenerate} implies $B$ is non-degenerate.
\end{proof}

With the previous Theorem we can extend Example \ref{example globalize al2ai} considerably.
Let $A$ be a *-algebra and $I$ an *-ideal of it.
Recall that an orthogonal complement for $I$ is a *-ideal $J\subset A$ such that $A=I\oplus J.$
In this situation $IJ=0$ because $IJ\subset I\cap J.$

Several remarks are in order.
Firstly, if $I$ is non degenerate the it has at most one orthogonal complement.
Secondly, if $J$ is an orthogonal complement for $I,$ then $A$ and $I$ are non degenerate $\Leftrightarrow$ $I$ and $J$ are non degenerate $\Leftrightarrow$ $A$ and $J$ are non degenerate.
Thirdly, $I$ is assimilative in $A$ every time $A$ and $I$ are non degenerate and $I$ has an orthogonal complement.
Finally, there are non degenerate assimilative ideals without an orthogonal complement.
As an example take the complex sequences vanishing except on a finite set, considered as a *-ideal of the complex sequences converging to zero.

\begin{corollary}
  Let $\al=\left(\{A_t\}_{t\in G},\{\al_t\}_{t\in G}\right)$ be a *-partial action of $G$ on $A.$
  If for all $s,t\in G$ the *-subalgebras $A_t+A_s$ and $A_t\cap A_s$ are non-degenerate and $A_t$ has an orthogonal complement in $A,$ then $\al$ has a non degenerate and minimal *-globalization.
\end{corollary}
\begin{proof}
  From the hyphoteses it follows that $A=A_e+A_e$ and $A_t=A_t+A_t$ are non degenerate.
  Besides $A_t$ has a unique orthogonal complement in $A,$ so it is assimilative in $A.$
  Then it suffices to verify that $\al$ satisfies  condition (2) from Theorem \ref{theorem enveloping star algebra}.
  
  Take $(t,a,b)\in G\times A\times A.$
  Since $A_\tmu$ has an orthogonal complement in $A,$ there exists $a'\in A_\tmu$ such that $ca=ca',$ for all $c\in A_\tmu.$
  With $u:=\al_t(a')b\in A_t$ we have, for all $c\in A_\tmu:$ $\al_t(c)u = \al_t(c)\al_t(a')b=\al_t(ca')b=\al_t(ca)b.$
\end{proof}

It may happen that a *-partial action has a minimal *-globalization with many ideals $A_t$ without an orthogonal complement.
To exhibit an example consider let $C_0(\R)$ be the *-algebra of continuous function from $\R$ to $\C$ vanishing at $\pm\infty.$
Let $\beta$ be the action of $\R$ on $C_0(\R)$ given by $\beta_t(f)(r) = f(r-t).$
As an ideal of $C_0(\R)$ take $A=C_0(0,+\infty)=\{f\in C_0(\R)\colon f|_{\R\setminus (0,+\infty)}\equiv 0\}$ and let $\alpha$ be the restriction of $\beta$ to $A.$
Then we can construct a non degenerate and minimal *-globalization of $\alpha$ using $\beta.$
Every ideal of $A$ is non degenerate because in $A$ the identity $aa^*=0$ implies $a=0,$ but the only *-ideals of $A$ with an orthogonal complement are the trivial ones.
For example, with $t>0$ the ideal $A_t=A\cap \be_t(A)=C_0(t,+\infty)$ does not have an orthogonal complement in $A.$

\subsubsection{About *-partial actions on commutative algebras}

We close this section with an adaptation of \cite[Proposition 2.1]{Ab03} to *-algebras.

\begin{proposition}
  Assume $\al$ is a *-partial action of $G$ on the commutative *-algebra $A$ and $\sigma=(B,\beta,I,\iota)$ a minimal *-globalization.
  If $A_t$ is non degenerate for all $t\in G$ then $B$ is commutative.
\end{proposition}
\begin{proof}
  Since $B=\spn\{\be_r(a)\colon r\in G,\ a\in I\},$ it suffices to prove that $\be_r(a)\be_s(b)=\be_s(b)\be_r(a),$ for all $a,b\in I$ and $r,s\in G.$
  Fix $a,b\in A$ and $r,s\in G.$
  With $t:=\rmu s$ we have $\be_r(a\be_t(b))=\be_r(a)\be_s(b)$ and $\be_r(\be_t(b)a)=\be_s(b)\be_r(a).$
  Then the proof will be completed once we show that $a\be_t(b)=\be_t(b)a.$
  
  Note that $I$ is commutative, $I_t:=\iota(A_t)$ non degenerate and $a\be_t(b),\be_t(b)a\in I\cap \be_t(I)=I_t.$
  Then $a\be_t(b)=\be_t(b)a$ if and only if $ca\be_t(b)=c\be_t(b)a,$ for all $c\in I_t.$
  Take $c\in I_t,$ using that $I$ is commutative and $a,b,c,\be_\tmu(ac),\be_\tmu(c),\be_t(b)a\in I,$ we obtain
  $ ca\be_t(b) = \be_t(\be_\tmu(ac)b)=\be_t(b\be_\tmu(ac))=\be_t(b)ac = c\be_t(b)a.$
\end{proof}

\section{Partial actions on C*-algebras}\label{section pa on cast algebras}

Assume $A$ is a C*-algebra and $G$ a topological group (we do not require it to be locally compact nor Hausdorff).

\begin{definition}
    The pair $\al=\left(\{A_t\}_{t\in G},\{\al_t\}_{t\in G}\right)$ is a C*-partial action (of $G$ on $A$) if:
    \begin{itemize}
     \item It is a *-partial action of $G$ on $A.$
     \item For all $t\in G,$ $A_t$ is closed in $A.$
     \item $\{A_t\}_{t\in G}$ is a continuous family, that is: for every open set $U\subset A$ the set $G_U:=\{t\in G\colon A_t\cap U\neq \emptyset\}$ is open in $G.$
     \item The function $\{(t,a)\in G\times A\colon a\in A_\tmu\}\to A,\ (t,a)\mapsto \al_t(a),$ is continuous.
    \end{itemize}
\end{definition}

Morphisms between C*-partial actions are just morphism between *-partial actions.
Recall from \cite[Example 2.1]{Ab03} that the restriction of a C*-partial action to a C*-ideal is a C*-partial action.

\begin{definition}
    A C*-globalization of the C*-partial action $\al,$ of $G$ on $A,$ is a 4-tuple $\Xi=(B,\beta,I,\iota)$ where:
    $B$ is a C*-algebra, $\beta$ is a C*-global action of $G$ on $B,$ $I$ is a C*-ideal of $B$ and $\iota\colon \al\to \be|_I$ is an isomorphism of *-partial actions.  
\end{definition}

Note we do not require $\iota$ to be a homeomorphism, this is automatic because every *-homomorphism between C*-algebras is contractive and has closed range.

In case $\al$ is a C*-partial action, it is also a *-partial action and every C*-globalization of $\al$ is a *-globalization of $\al.$
Although not every *-globalization of $\al$ is a C*-globalization, the existence of a *-globalization implies the existe of a C*-globalization.
We will prove this claim in two steps: first we show the group's topology is irrelevant to globalize C*-partial actions, then we construct a C*-globalization from a *-globalization.

In the next statement $G^\dis$ denotes the group $G$ with the discrete topology.

\begin{lemma}
    Let $G$ be a topological group, $\be$ a C*-global action of $G^\dis$ on $B$ and $I$ an ideal of $B$ such that $[I]:=\spn\{ \be_t(I)\colon t\in G \}$ is dense in $B.$
    Set $\al:=\be|_I,$ which is a C*-partial action of $G^\dis$ on $I.$
    Then $\be$ is a C*-partial action of $G$ if and only if $\al$ is a C*-partial action of $G.$
\end{lemma}
\begin{proof}
    The direct implication is \cite[Example 2.1]{Ab03}.
    For the converse note that (as each $\be_r$ is an isometry) $G\times B\to B,\ (t,b)\mapsto \be_t(b),$ is continuous if and only if for every $b\in B$ the function $ev_b\colon G\to B,\ t\mapsto \be_t(b),$ is continuous at $e.$
    Moreover, $U:=\{b\in B\colon ev_b \mbox{ is continuous at } e\}$ is a closed $\be$-invariant subspace of $B,$ thus all we need to show is that $I\subset U.$
    
    Fix $a\in I.$
    From \cite[Lemma 2.1]{Ab03} (using the family of ideals $\{\be_t(I)\}_{t\in G}$) we get, for all $t\in G,$ that 
    \begin{align*}
	\|\be_t(a)-a\|
	    &= \sup\{ \|(\be_t(a)-a)\be_r(b)\|\colon r\in G,\ b\in I,\ \|b\|\leq 1 \}\\
	    & = \sup\{ \|(\be_{\rmu t}(a)-\be_{\rmu}(a))b\|\colon r\in G,\ b\in I,\ \|b\|\leq 1 \}\\
	    & = \sup\{ \|(\be_{r t}(a)-\be_{r}(a))b\|\colon r\in G,\ b\in I,\ \|b\|\leq 1 \}.
    \end{align*}

    Fix $r\in G$ and $b\in I$ with $\|b\|\leq 1.$
    Note that $(\be_{r t}(a)-\be_{r}(a))b = \be_{r t}(a)b-\be_{r}(a)b \in \be_{rt}(I)I + \be_r(I)I = I_{rt}+ I_r.$ 
    Thus, if $\{v_\lambda\}_{\lambda\in \Lambda}$ is an approximate unit of $I_{rt}+ I_r,$ we have $\|(\be_{r t}(a)-\be_{r}(a))b\|=\lim_\lambda\|(\be_{r t}(a)-\be_{r}(a))bv_\lambda\|.$
    
    Given $\lambda\in \Lambda,$ $v_\lambda\in (I_{rt}+ I_r)^+$ and there exists $c\in I_{rt}^+$ and $d\in I_r^+$ such that $v_\lambda = c+d.$
    This implies $\|c\|,\|d\|\leq \|v_\lambda\|\leq 1.$
    One one hand, if $\{w_\mu\}_{\mu\in M}$ is an approximate unit of $I_t,$ we have
    \begin{align*}
	\|(\be_{r t}(a)-\be_{r}(a))bc\|
	    & = \| a\be_{\tmu\rmu}(bc) - \be_{\tmu}(a)\be_{\tmu\rmu}(bc) \|\\
	    & =\lim_\mu \| a\be_{\tmu\rmu}(bc) - \be_\tmu(w_\mu)\be_{\tmu}(a)\be_{\tmu\rmu}(bc) \|\\
	    & \leq \limsup_\mu \| a- \be_\tmu(w_\mu a)\|
		=\limsup_\mu \| a- \al_\tmu(w_\mu a)\|.
    \end{align*}
    
    To construct a particular approximate unit consider $M_t:=\{ c\in I_t^+\colon \|a\|<1\}$ with it's natural order, then $\{\mu\}_{\mu\in M_t}$ is an approximate unit of $I_t$ and
    $$ \|(\be_{r t}(a)-\be_{r}(a))bc\| \leq \lim_{\mu\in M_t} \sup\{ \| a- \al_\tmu(\nu a)\|\colon \nu\in M_t,\nu \geq \mu\}=: C(t).$$
    
    With $s=rt$ by symmetry we get
    \begin{align*}
	\|(\be_{r t}(a)-\be_{r}(a))bd\| = \|(\be_{s\tmu }(a) - \be_{s}(a) )bd\| \leq C(\tmu).
    \end{align*}
    Putting all together we obtain $\|\be_t(a)-a\|\leq C(t) + C(\tmu).$
    All we need to show is that $\lim_{t\to e} C(t)=0.$
    
    Take $\vep>0.$
    As $\al$ is a partial action of $G$ on $I$ there are neighbourhoods $V\subset G$ and $W\subset B$ of $e$ and $a,$ respectively, such that for $s\in V$ and $b\in I_\smu\cap W$ we have $\| \al_s(b) - a \|<\vep/2.$
    Take $\delta>0$ such that $B(a,\delta)\subset W.$
    Then $U:=\{r\in G\colon B(a,\delta/2)\cap I_r\neq \emptyset \}$ is an open set containing $e$ because $\{I_t\}_{t\in G}$ is a continuous family.
    For $r\in U\cap V^{-1}$ there exists $b\in I_r\cap B(a,\delta/2),$ so $\lim_{\mu\in M_r} \|a-\mu a\|=dist(a,I_r)\leq \delta/2.$
    Take $\mu_r\in M_r$ such that $\| a - \nu a \|<\delta$ for all $\nu\in M_r$ with $\nu\geq\mu_r.$
    Then for $r\in U\cap V^{-1}$ we have $C(r)<\vep$ because for all $\nu\in M_r$ with $\nu \geq \mu_r$ the inequality $\|a - \al_\rmu(\nu a)\|\leq \vep/2$ holds and implies
    $$ \sup\{ \|a - \al_\rmu(\nu' a)\|\colon \nu'\in M_r,\ \nu'\geq \nu \} \leq \vep/2 <\vep. $$
\end{proof}

\begin{corollary}\label{corollary topology irrelevant}
    Let $\al=\left(\{A_t\}_{t\in G},\{\al_t\}_{t\in G}\right)$ be a C*-partial action and use the expression $\al^\dis$ to denote $\al$ as a C*-partial action of $G^\dis.$
    Then $\al$ has a C*-globalization if and only if $\al^\dis$ has a C*-globalization.
\end{corollary}
\begin{proof}
    As mentioned before the direct implication in immediate.
    For the converse assume $\al^\dis$ has a C*-globalization $\Xi=(B,\be,I,\iota).$
    Without loss of generality we may assume $\Xi$ is minimal (an enveloping action in the sense of F. Abadie \cite{Ab03}) this is $\overline{[I]}=B.$
    Note $\be|_I$ is a C*-partial action of $G$ because it is isomorphic (as a *-partial action) to $\al.$
    Then the previous lemma implies $\be$ is a C*-global action of $G.$
    Thus $\Xi$ is a C*-globalization of $\al.$
\end{proof}

As promised before we now construct a C*-globalization from a *-globalization.

\begin{theorem}\label{main theorem c star algebras}
    Let $\al=\left(\{A_t\}_{t\in G},\{\al_t\}_{t\in G}\right)$ be a C*-partial action.
    Then the following are equivalent:
    \begin{enumerate}[(1)]
     \item $\al$ has a C*-globalization.
     \item $\al$ has a *-globalization.
     \item For all $(t,a,b)\in G\times A\times A$ there exists $u\in A_t$ such that, for all $c\in A_\tmu,$ $\al_t(c)u=\al_t(ca)b.$
    \end{enumerate}
\end{theorem}
\begin{proof}
    We already know (1) implies (2).    
    On C*-algebras the condition $a^*a=0$ implies $a=0,$ thus every C*-ideal of a C*-algebra is non-degenerate.
    The existence of approximate units implies every closed ideal is assimilative.
    Besides the sum\footnote{Without closure.} of two C*-ideals of a C*-algebra is, again, a C*-ideal.
    Then (2) and (3) are equivalent by Theorem \ref{theorem enveloping star algebra}.
    
    To show (2) implies (1), without loss of generality, we assume $G$ is discrete (Corollary \ref{corollary topology irrelevant}).
    From the previous paragraph and Proposition \ref{proposition of uniqueness} we know $\al$ has a non degenerate and enveloping *-globalization $\Xi=(B,\beta,I,\iota).$    
    Note $I$ is a C*-algebra because it is isomorphic (as a *-algebra) to a C*-algebra.
    Moreover, $\be|_I$ is a C*-partial action isomorphic to $\al.$
    Then it suffices to show that $\be|_I$ has a C*-globalization and me may think $I=A,$ $\iota$ is inclusion of $A$ in $B$ and $\al=\be|_I.$
    
    Let $\pi$ be the canonical morphism associated to $\Xi$ and $\Theta$ the canonical action of $G$ on $M(A^G).$
    From Proposition \ref{proposition of uniqueness} we know $\pi$ is injective and, as we think $A\subset B,$ $\pi(b)f|_r = \be_t(b)f|_r.$
    
    The set of bounded functions from $G$ to $A,$ $A^G_b,$ is a C*-algebra with the *-algebra structure inherited from $A^G$ and the supremum norm.
    Define 
    $$C:=\{T\in M(A^G)\colon T(A^G_b)\cup T^*(A^G_b)\subset A^G_b\}.$$
    
    Note $C$ is $\Theta$ invariant and that to show $\pi(B) \subset C$ it suffices to prove that $\pi(A)\subset C.$
    This last inclusion holds because given $a\in A,$ $f\in A^G_b$ and $r\in G$ and an approximate unit of $A_\rmu,$ $\{v_\lambda\}_{\lambda\in \Lambda},$ we have $\be_r(a)f|_r\in A_r$ and
    $$\|\pi(a)f|_r\|=\|\be_r(a)f|_r\|=\lim_\lambda \|\al_r(v_\lambda)\be_r( a) f|_r\|=\lim_\lambda \| \al_r(v_\lambda a) f|_r \|\leq \|a\|\|f\|.$$
    
    Consider $M(A^G_b)$ as a C*-algebra and let $\rho\colon \pi(B)\to M(A^G_b)$ be defined as $\rho(T)f = Tf.$
    To show $\rho$ is injective it suffices to show that $\rho\circ \pi $ is injective.
    Take $b\in B$ such that $\rho\circ\pi(b)=0.$
    There are $a_1,\ldots,a_n\in A$ and $t_1,\ldots,t_n\in G$ such that $b=\sum_{j=1}^n \be_{t_j}(a_j).$
    Given $r\in G$ and $c\in A$ we have $ b\be_r(c) = \be_r(\be_\rmu(b)c) = \be_r( \rho\circ\pi(b)\delta_\rmu^c|_\rmu )=0.$
    Thus $bB=\{0\}$ and this implies $b=0.$
    
    Given $t\in G$ set $\psi_t\colon A^G_b\to A^G_b$ as $\psi_t(f)|_r = f|_{rt}$ and $\Psi_t\colon M(A^G_b)\to M(A^G_b)$ as $\Psi_t(T)=\psi_t\circ T\circ \psi_\tmu.$
    Then $\Psi$ is a C*-global action of $G$ on $M(A^G_b);$ $\rho\circ \pi\colon \be\to \Psi$ is a morphism of C*-partial actions and the closure of $\rho\circ \pi(B),$ $D,$ is a $\Psi$-invariant C*-subalgebra of $M(A^G_b).$
    Name $\gamma$ the restriction of $\Psi$ to $D.$
    Note $J:=\rho(\pi(A))$ is a C*-ideal of $D$ because $\rho\circ \pi|_A$ has closed range and $J$ is an ideal of $\rho\circ \pi(B).$
    From Remark \ref{remark double restriction} we get $\gamma|_J=\Psi|_D|_{\rho\circ \pi(A)}=\Psi|_{\rho\circ \pi(A)}=\Psi|_{\rho\circ\pi(B)}|_{\rho\circ\pi(A)}.$
    Besides, $\rho\colon \Theta|_{\pi(B)}\to \Psi|_{\rho\circ \pi(B)}$ is an isomorphism.
    Then $\rho\circ \pi|_A\colon \al \to \Psi|_J$ is an isomorphism and, consequently, $(D,\gamma,J,\rho\circ \pi|_A)$ is a C*-globalization of $\al.$
\end{proof}

Condition (3) asserts the existence of a certain element $u$ for each $(t,a,b)\in G\times A\times A.$
That element is the limit of the net given in (2) below.

\begin{proposition}\label{prop equivalence globalization}
     Let $\al=\left(\{A_t\}_{t\in G},\{\al_t\}_{t\in G}\right)$ a C*-partial action.
     Then the following are equivalent.
     \begin{enumerate}[(1)]
      \item $\al$ has a C*-globalization.
      \item There exists $U,V\subset A$ such that: (i) $\spncl AU=\spncl VA=A$ and (ii) for all $(t,a,b)\in G\times U\times V$ there exists an approximate unit of $A_\tmu,$ $\{v_\lambda\}_{\lambda\in \Lambda},$ such that $\{\al_t(v_\lambda a)b\}_{\lambda\in \Lambda}$ is a Cauchy net (or converges).
     \end{enumerate}
\end{proposition}
\begin{proof}
    In case (1) holds set $U=V=A$ and take $(t,a,b)\in G\times A\times A.$
    Let $u$ be the element given in (2) of the previous theorem and $\{v_\lambda\}_{\lambda\in \Lambda}$ and approximate unit of $A_\tmu.$
    Then $\{\al_t(v_\lambda)\}_{\lambda\in \Lambda}$ is an approximate unit of $A_t$ and $\{\al_t(v_\lambda)u\}_{\lambda\in \Lambda}=\{\al_t(v_\lambda a)b\}_{\lambda\in \Lambda}$ converges to $u.$
    
    Now assume (2) is true.
    Take $(t,a,b)\in G\times A\times A$ for which there exists an approximate unit, $\{v_\lambda\}_{\lambda\in \Lambda},$ such that $\{\al_t(v_\lambda a)b\}_{\lambda\in \Lambda}$ converges to $u.$
    We claim that given any other approximate unit of $A_\tmu,$ $\{w_\mu\}_{\mu\in M},$ the net $\{\al_t(w_\mu a)b\}_{\mu\in M}=\{\al_t(w_\mu)u\}_{\mu\in M}$ converges to $u.$
    Indeed
    $$ \lim_\lambda \al_t(v_\lambda a)b = \lim_\mu \al_t(w_\mu)\lim_\lambda \al_t(v_\lambda a)b =\lim_\mu \lim_\lambda \al_t(w_\mu v_\lambda a)b =\lim_\mu \al_t(w_\mu a)b.$$
    Then $u$ is determined by $(t,a,b),$ so we denote it $u_t(a,b).$
    Moreover, if $c\in A_\tmu$ then
    $$\al_t(c)u_t(a,b) =\lim_\lambda \al_t(c)\al_t(v_\lambda a)b = \lim_\lambda \al_t(cv_\lambda a)b =\al_t(ca)b. $$
    
    Take $a\in U,$ $b\in V$ and $c,d\in A$ and any approximate unit of $A_\tmu,$ $\{v_\lambda\}_{\lambda\in \Lambda}.$
    Then $\{\al_t(v_\lambda ca)bd\}_{\lambda\in \Lambda}$ converges to $\al_t( c\al_\tmu( u_t(a,b)d ) )$ because, for all $\lambda\in \Lambda,$
    $$ \al_t(v_\lambda ca)bd = \al_t(v_\lambda c)u_t(a,b)d = \al_t(v_\lambda c \al_\tmu(u_t(a,b)d ) ) $$
    and $\lim_\lambda \al_t(v_\lambda c \al_\tmu(u_t(a,b)d ) )= \al_t(c \al_\tmu(u_t(a,b)d ) ).$
    
    The conclusions of the previous paragraphs imply that (2) holds if we replace $U$ for $AU$ and $V$ for $VA.$
    From now on we assume $AU=U$ and $VA=V.$
    
    Given $a,b\in U,$ $c\in V$ and $\lambda\in \C$ note that $u_t(a,c)+\lambda u_t(b,c)\in A_t$ and for all $d\in A_\tmu:$
    $ \al_t(d)(u_t(a,b)+\lambda u_t(b,c))=\al_t(d(a+\lambda b))c.$
    In case $a\in U$ and $b,c\in V:$ $\al_t(d)(u_t(a,b)+\lambda u_t(a,c))=\al_t(da)(b+\lambda c).$
    Then, for all $t\in G,$ there exists a unique bilinear function $u_t\colon \spn U\times \spn V\to A_t$ such that for all $a\in \spn U,\ b\in \spn V$ and $c\in A_\tmu: $  $\al_t(c)u_t(a,b)=\al_t(ca)b.$
    Note that $\|u_t(a,b)\|^2=\|\al_t( \al_\tmu(u_t(a,b)^*) a)b\|\leq \|u_t(a,b)\|\|a\|\|b\|,$ then $\|u_t(a,b)\|\leq \|a\|\|b\|.$
    So there is a unique continuous bilinear function $v_t\colon A\times A\to A_t$ extending $u_t,$ for every $t\in G.$
    
    Given $(t,a,b)\in G\times A\times A$ and $c\in A_\tmu$ take sequences $\{a_n\}_{n\in \N}\subset \spn U$ and $\{b_n\}_{n\in \N}\subset \spn V$ converging to $a$ and $b,$ respectively.
    Then 
    $$\al_t(c)v_t(a,b)=\lim_{n}\al_t(c)v_t(a_n,b_n)=\lim_n \al_t(ca_n)b_n = \al_t(ca)b.$$
    This shows (2) implies (3) of the previous theorem, so $\al$ admits a C*-globalization.
\end{proof}

Now we use our criterion to give some conditions were a C*-partial action can be globalized.

On a unital C*-algebra, $A,$ for every unital ideal $I$ there exists a unique ideal $J$ such that $A=I\oplus J$ and $IJ=\{0\}.$
In general, even for non unital $A,$ we say a C*-ideal $I$ has an orthogonal complement if there exists a C*-ideal $J$ such that $A=I\oplus J$ and $IJ=\{0\}.$

\begin{corollary}
    Let $\al=\left(\{A_t\}_{t\in G},\{\al_t\}_{t\in G}\right)$ be a C*-partial action.
    If, for all $t\in G,$ $A_t$ has an orthogonal complement then $\al$ has a C*-globalization.
\end{corollary}
\begin{proof}
    Given $(t,a,b)\in G\times A\times A$ let $J$ be the orthogonal complement of $A_\tmu.$
    Then there are $a_1\in A_\tmu$ and $a_2\in J$ such that $a=a_1+a_2.$
    Note for all $c\in A_\tmu,$ $ca =ca_1.$
    Thus $u:=\al_t(a_1)b\in A_t$ and $\al_t(c)u=\al_t(ca_1)b = \al_t(ca)b.$
\end{proof}

The next corollary is a version for C*-algebras of \cite[Theorem 4.5]{ExDo05}.
The proof is omitted because it follows form the cited theorem and the previous corollary. 

\begin{corollary}
    Let $\al=\left(\{A_t\}_{t\in G},\{\al_t\}_{t\in G}\right)$ a C*-partial action with $A$ unital.
    Then $\al$ has a C*-globalization if and only if $A_t$ is a unital algebra, for all $t\in G.$
\end{corollary}

There are certain *-homomorphisms that force the existence of C*-enveloping actions.
For partial actions on locally compact and Hausdorff spaces the next result can be showed using \cite[Proposition 2.1]{Ab03} instead of Theorem \ref{main theorem c star algebras} (in that case the homomorphism $\phi$ defines a continuous function between the spectra of the algebras).

\begin{corollary}
  Let $\al=\left(\{A_t\}_{t\in G},\{\al_t\}_{t\in G}\right)$ and $\be=\left(\{B_t\}_{t\in G},\{\be_t\}_{t\in G}\right)$ be C*-partial actions.
  If $\al$ has a C*-globalization and there exists a *-algebra's homomorphism $\phi\colon A\to M(B)$ such that:
  \begin{enumerate}
    \item for all $t\in G,$ $\spncl \{\phi(a)b\colon a\in A_t,\ b\in B\} = B_t$ and
    \item for all $t\in G,$ $a\in A_\tmu$ and $b\in B_\tmu,$ $\phi(\al_t(a))\be_t(b)=\be_t(\phi(a)b);$
  \end{enumerate}
  then $\be$ has a C*-globalization.
\end{corollary}
\begin{proof}
  By Cohen-Hewitt's Theorem $B_t =\{\phi(a)b\colon a\in A_t,\ b\in B\},$ for all $t\in G.$
  The same theorem implies $A_t = \{ab\colon a,b\in A_t\},$ then $B_t=\{\phi(a)b\colon a\in A_t,\ b\in B_t\}.$
  
  Given $(t,b_1,b_2)\in G\times B\times B$ take $a_1,a_2\in A$ and $c_1,c_2\in B$ such that $b_1^* =\phi(a_1^*)c_1^* $ and $b_2=\phi(a_2)c_2.$
  Now let $u$ be the element given for $(t,a_1,a_2)$ by condition (3) of Theorem \ref{main theorem c star algebras}.
  Then $c_1\phi(u)c_2\in B_t$ because $u\in A_t.$
  
  Given $d\in B_\tmu$ choose $d'\in B_\tmu$ and $e\in A_\tmu$ such that $\be_\tmu(c_1^*\be_t(d^*)) = \phi(e^*)d'^*.$
  Then
  \begin{align*}
    \be_t(d) c_1\phi(u)c_2
      & = \be_t(d')\phi(\al_t(e)u)c_2
          = \be_t(d')\phi( \al_t(ea_1) a_2 ) c_2\\
      & = \be_t(d')\phi(\al_t(ea_1)) \phi(a_2)c_2
          = \be_t( d'\phi(e)\phi(a_1) ) b_2 =\be_t(db_1)b_2 .
  \end{align*}
  The rest follows directly from the previous theorem.
\end{proof}

To close this section we give an example to show condition (1) from the previous statement can not we weakened to (1') $B_t =\spncl \{\phi(a)b\colon a\in A_t,\ b\in B_t\},$ for all $t\in G.$

Let $G= \Zb_2 = \{1,-1\}$ (with multiplicative notation) $A=B=C[0,1].$
As $\al$ consider the trivial global action of $G$ on $A.$
To give a partial action of $G$ on $B,$ $\be,$ we just need to specify $\be_{-1},$ which will be the identity on $C_0([0,1)).$
For $\phi\colon C[0,1]\to M(C[0,1])=C[0,1]$ just take the identity. 
Note (1') and (2) are satisfied, $\al$ has a C*-globalization but $\be$ does not.

\section{Partial action on equivalence bimodules}\label{section pa on equiv bimodules}
In \cite{Ab03} F. Abadie defines Morita equivalence of C*-partial actions using partial actions on positive C*-trings.
Recall that a positive C*-trings are exactly equivalence bimodules \cite{Ab03,Zl83,Rf82Morita}.
Here we adopt the terminology of equivalence bimodules.

In this last section we give a necessary and sufficient condition for the existence of a globalization of a partial action on an equivalence bimodule and, as a consequence, we obtain the uniqueness of Morita enveloping actions (as was shown in \cite{Ab03}).

We adopt the terminology of \cite{Rf82Morita} and agree that ``${}_A\X_B$ is an equivalence bimodule'' means ``$\X$ is an $A$-$B$-equivalence bimodule''.
Take equivalence bimodules ${}_A\X_B$ and a ${}_C\Y_D.$
A function $\phi\colon \X\to \Y$ is an Hb-homomorphism if it is linear and for all $x,y,z\in \X,$ $\phi(x\la y,z\ra_B)=\phi(x)\la \phi(y),\phi(z)\ra_D.$
Such functions are contractive and, with the previous notation, there exist unique *-homomorphisms ${}^l\phi\colon A\to C$ and $\phi^r\colon B\to D$ such that, for all $x,y\in \X,$ $\phi^l({}_A\la x,y\ra)={}_C\la \phi(x),\phi(y)\ra$ and $\phi^r(\la x,y\ra_B)=\la \phi(x),\phi(y)\ra_D.$
The proof of these facts can be found in \cite{Ab03}.

Given a C*-ideal $I$ of $A,$ by Cohen-Hewitt $I\X:=\{ax\colon a\in I,\ x\in \X\}$ is a closed submodule of $\X.$
We denote ${}_IB$ the C*-ideal of $B$ induced by $I$ (or $I\X$) through $ \X,$ that is ${}_IB=\spncl \{\la u,v \ra_B\colon u,v\in I\X \}.$
In a similar way we define, for a C*-ideal $J$ of $B,$ $\X J$ and $A_J.$
Recall that $A_{{}_IB}=I$ and $J={}_{A_J}B.$

For a closed subspace $ \Z$ of $\X$ the following are equivalent: (i) there exists a C*-ideal $I$ of $A$ such that $\Z = I\X,$ (ii) there exists a C*-ideal $J$ of $B$ such that $\Z = \X J,$ (iii) $\X \la \Z,\X\ra_B\subset \Z$ and (iv) ${}_A\la \X,\Z\ra \X\subset \Z.$
If these conditions are satisfied then $A_\Z:=\spncl {}_A\la \Z,\Z\ra$ and ${}_{\Z}B$ are C*-ideals, $\Z=A_{\Z}\X=\X{}_{\Z}B$ and we say $\Z$ is an ideal of $\X.$
Every ideal $\Z$ of $\X$ is an $A_\Z-{}_\Z B$-equivalence bimodule.

\begin{definition}
  The pair $\gamma = \left(\{\gamma_t\}_{t\in G},\{\X_t\}_{t\in G}\right)$ is a \textit{Hb-partial action} if:
  \begin{itemize}
   \item $\X$ is an ($A$-$B$-)equivalence bimodule and $G$ a topological group.
   \item $\gamma$ is a set theoretic partial action of $G$ on $\X.$
   \item $\{\X_t\}_{t\in G}$ is a continuous family of ideals of $\X.$
   \item For all $t\in G,$ $\gamma_t\colon \X_\tmu\to \X_t$ is an Hb-homomorphism.
   \item The function $\{(t,x)\in G\times \X\colon x\in \X_\tmu\}\to \X,\ (t,x)\mapsto \gamma_t(x),$ is continuous.
  \end{itemize}
\end{definition}

Assume $\delta$ is a Hb-partial action $G$ on ${}_C\Y_D.$
We say $\phi\colon \gamma\to \delta$ is a \textit{Hb-morphism} if it is a Hb-homomorphism from $\X$ to $\Y$ which is also a morphism of set theoretic partial actions.

From \cite{Ab03} we know there are unique C*-partial actions, $\al = \left(\{\al_t\}_{t\in G},\{A_t\}_{t\in G}\right)$ and $\be = \left(\{\be_t\}_{t\in G},\{B_t\}_{t\in G}\right),$ such that
\begin{itemize}
 \item For all $t\in G,$ $A_t = A_{\X_t}$ and $B_t={}_{\X_t}B.$
 \item For all $t \in G$ and $x,y\in \X_\tmu:$
 $$ \al_t({}_A\la x,y\ra)={}_A\la \gamma_t(x),\gamma_t(y)\ra \mbox{ and }\be_t(\la x,y\ra_B)=\la \gamma_t(x),\gamma_t(y)\ra_B.$$
\end{itemize}

We will call $\al$ the left side of $\gamma$ and $\be$ the right side of $\gamma$ and will denote them ${}^l\gamma$ and $\gamma^r,$ respectively.

\begin{example}
  Every C*-partial action, $\al,$ on a C*-algebra $A$ is a Hb-partial action on ${}_AA_A.$
  Besides $\al={}^l\al=\al^r.$
\end{example}

\begin{example}\label{example of restriction}
  Given a Hb-global action of $G$ on ${}_A\X_B,$ $\gamma,$ and an ideal $\Y$ of $\X;$ the restriction $\gamma|_\Y$ is a Hb-partial action on the $A_\Y$-${}_\Y B$-equivalence bimodule $\Y.$
  In this case  ${}^l(\gamma|_\Y)=({}^l\gamma)|_{A_\Y}$ and $(\gamma|_\Y)^r = \gamma^r|_{{}_\Y B}.$
\end{example}

With restrictions of global actions, on one hand, and isomorphisms on the other we define globalizations.

\begin{definition}
  Let $\gamma$ be an Hb-partial action of $G$ on $\X.$
  A globalization of $\gamma$ is a $4-$tuple $\Xi=(\Y,\delta,\Z,\iota)$ such that:
  $\Y$ is an equivalence bimodule, $\delta$ is an Hb-global action of $G$ on $\Y,$ $\Z$ is an ideal of $\Y$ and $\iota\colon \gamma\to \delta|_{\Z}$ is an isomorphism of Hb-partial actions.
  In case $[\Z]:=\spn \{\delta_t(\Z)\colon t\in G\}$ is dense in $\Y$ we say $\Xi$ is a minimal globalization.
\end{definition}

The nexus between Hb-partial actions and C*-partial actions is the linking partial action \cite{Ab03}.
To describe this action we start with an Hb-partial action of a group $G$ on ${}_A\X_B,$ $\gamma,$ and set $\al:={}^l\gamma$ and $\be:=\gamma^r.$
The linking algebra of $\X$ is the algebra of generalized compact operators of the $A$-Hilbert module $\X\oplus A,$ $\Lb(\X)=\Kb(\X\oplus A).$
In matrix representation $ \Lb(\X) =\left( \begin{smallmatrix} A & \X \\ \widetilde{\X} & B \end{smallmatrix} \right), $ with $\widetilde{\X}$ the $B$-$A$-equivalence bimodule adjoint to $\X.$

The linking partial action of $\gamma,$ $\Lb(\gamma)=\left(\{\Lb(\gamma)_t\}_{t\in G},\{ \Lb(\X)_t \}_{t\in G} \right),$ is the unique C*-partial action such that, for all $t\in,$ $\Lb(\X)_t=\Lb(\X_t)$ and 
$$\Lb(\gamma)_t\left( \begin{smallmatrix} a & x\\ \widetilde{y} & b \end{smallmatrix} \right) = \left( \begin{smallmatrix} \al_t(a) & \gamma_t(x)\\ \widetilde{\gamma_t(y)} & \be_t(b) \end{smallmatrix} \right)\ \forall\ x,y\in X_\tmu,\ a\in A_\tmu\mbox{ and }b\in B_\tmu.$$

In case $\delta$ is an Hb-partial action of $G$ on $\Y$ and $\pi\colon \gamma\to \delta$ is an isomorphism, the morphism $\Lb(\pi)\colon \Lb(\gamma)\to \Lb(\delta)$ defined as
$$ \Lb(\pi)\left( \begin{smallmatrix} a & x\\ \widetilde{y} & b \end{smallmatrix} \right) = \left( \begin{smallmatrix} {}^l\pi(a) & \pi(x)\\ \widetilde{\pi(y)} & \pi^r(b) \end{smallmatrix} \right),$$
is an isomorfphism with inverse $\Lb(\pi^{-1}).$

\begin{proposition}\label{prop globalization of pa on equiv bimod}
    If $\Xi=({}_C\Y_D,\delta,\Z,\iota)$ is an enveloping globalization of $\gamma$ then
    \begin{itemize}
     \item $(\Lb(\Y),\Lb(\delta),\Lb(\Z),\Lb(\iota))$ is a C*-enveloping globalization of $\Lb(\gamma).$
     \item $(C,{}^l\delta,C_{\Z},{}^l\iota)$ is a C*-enveloping globalization of ${}^l\gamma.$
     \item $(D,\delta^r,{}_{\Z}D,\iota^r)$ is a C*-enveloping globalization of $\gamma^r.$
    \end{itemize}
\end{proposition}
\begin{proof}
  It is straightforward and is left to the reader.
\end{proof}

In the same way equivalence bimodules are constructed form C*-algebras and projections, Hb-partial actions are constructed form C*-partial actions and equivariant projections.

\begin{theorem}
 Let $\al$ be a C*-partial action of $G$ on $A$ and $p\in M(A)$ a projection such that $\spncl ApA =A$ and $\al_t(pa)=p\al_t(a),$ for all $t\in G$ and $a\in A_\tmu.$ 
 If $\X:=(1-p)Ap,$ $C:=(1-p)A(1-p)$ and $D:=pAp$ then
 \begin{enumerate}[(1)]
  \item $\X$ is a $C$-$D$-equivalence bimodule,
  \item $\X$ is $\al$-invariant,
  \item the restriction of $\al $ to $\X,$ $\gamma,$ is an Hb-partial action.
  \item $C$ and $D$ are $\al$-invariant, ${}^l\gamma=\al|_C$ and $\gamma^r=\al|_D.$
 \end{enumerate}
 
 Moreover, $\al$ is isomorphic (as a C*-partial action) to $\Lb(\gamma)$ and $\al$ has a C*-globalization if and only if $\gamma$ has an Hb-globalization.
\end{theorem}
\begin{proof}
  The proof of claims (1)-(4) are left to the reader.
  Besides, the usual identification of $A$ with $\Lb(\X)$ is an isomorphism of C*-partial actions between $\al$ and $\Lb(\gamma).$  
  Then, in case $\gamma$ has an Hb-globalization, $\Lb(\gamma)$ has a C*-globalization and this implies $\al$ has a C*-globalization.
  
  To prove the converse assume $\al$ has a C*-globalization.
  To prove $\gamma$ has a globalization we can assume, without loss of generality, that there exists a C*-partial action, $\be$ of $G$ on $B,$ such that $A$ is an ideal of $B,$ $\al=\be|_A$ and $\overline{[A]}=B.$
  
  The key claim is that there exists a unique projection $\overline{p}\in M(B)$ such that $\be_t(\overline{p}b)=\overline{p}\be_t(b)$ and $\overline{p}a=pa,$ for all $t\in G,$ $b \in B$ and $a\in A.$  
  To prove the claim it suffices to show that for all $t_1,\ldots,t_n\in G$ and $a_1,\ldots,a_n\in A$
  \begin{equation}
    \|\sum_{j=1}^n \be_{t_j}(pa_j) \|\leq \| \sum_{j=1}^n \be_{t_j}(a_j) \|.
  \end{equation}
  
  From \cite[Lemma 2.1]{Ab03} we conclude that it is enough to show that $$\|\sum_{j=1}^n \be_{t_j}(pa_j)\be_r(b) \|\leq \| \sum_{j=1}^n \be_{t_j}(a_j) \|,$$ for all $r\in G$ and $b\in A$ with $\|b\|<1.$
  
  Take $r\in G$ and $b\in B$ as before.
  Note $pa_j\be_{t_j^{-1} r}\in A_{t_j^{-1} r}$ so that
  \begin{align*}
    \|\sum_{j=1}^n \be_{t_j}(pa_j)\be_r(b) \|
      & = \|\sum_{j=1}^n \be_{t_j}(pa_j\be_{t_j^{-1} r}(b)) \|
	= \|\sum_{j=1}^n \be_{\rmu t_j}(pa_j\be_{t_j^{-1} r}(b)) \|\\
      & = \|\sum_{j=1}^n \al_{\rmu t_j}(pa_j\be_{t_j^{-1} r}(b)) \|
	= \|\sum_{j=1}^n p\al_{\rmu t_j}(a_j\be_{t_j^{-1} r}(b)) \|\\
      & \leq \|\sum_{j=1}^n \al_{\rmu t_j}(a_j\be_{t_j^{-1} r}(b))\|
	= \|\sum_{j=1}^n \be_{\rmu t_j}(a_j\be_{t_j^{-1} r}(b))\|\\
      & \leq  \|\sum_{j=1}^n \be_{t_j}(a_j)\|.
  \end{align*}
  
  Set $\Y:=(1-\overline{p})A\overline{p}\subset B,$ note that $\X\subset \Y$ and define $\iota\colon \X\to \Y$ as the canonical inclusion.
  Remark \ref{remark double restriction} implies $\be|_\Y|_\X=\be|_\X=\be|_A|_\X=\al|_\X=\gamma,$ then $(\Y,\be|_\Y,\X,\iota)$ is a globalization of $\gamma.$
\end{proof}

\begin{corollary}
    An Hb-partial action has a globalization if and only if it's linking partial action has a C*-globalization.
\end{corollary}
\begin{proof}
    Let $\gamma$ be a partial action of $G$ on the $A$-$B$-equivalence bimodule $\X$.
    The thesis follows easily from the previous theorem with $\al = \Lb(\gamma)$ and $p=\left(\begin{smallmatrix} 1&0\\0 & 0 \end{smallmatrix}
  \right)$ because $\Lb(\gamma)|_{p\Lb(\X)(1-p)}$ is isomorphic to $\gamma.$
\end{proof}

As a consequence of the Corollary we get that the group's topology does not affects the existence of Hb-globalizations because it does nor affects the existence of C*-globalizations.
This conclusion can be equally derived from our last result, which also implies that a C*-partial action has a C*-globalization if and only if it has a Hb-globalization.

\begin{theorem}
    An Hb-partial action has an Hb-globalization if and only if it's left and right sides have C*-globalizations.
\end{theorem}
\begin{proof}
    The direct implication follows from Proposition \ref{prop globalization of pa on equiv bimod}.
    For the converse assume $\gamma$ is an Hb-partial action of the topological group $G$ on ${}_A \X_B.$
    Assume $\al:={}^l\gamma$ and $\be:=\gamma^r$ have C*-globalizations.
    Given $t\in G$ set $A_t:=A_{\X_t}$ and $B_t:={}_{\X_t}B.$
    
    It suffices to show $\Lb(\gamma)$ has a C*-globalization, for which we use Proposition \ref{prop equivalence globalization} with $U=V=\{ \left(\begin{smallmatrix}
                   0 & x\\ \widetilde{y} & 0
                  \end{smallmatrix} \right)\colon x,y\in \X\}.$
    Take $(t,\xi,\eta)\in G\times U\times V.$
    Let $\{v_\lambda\}_{\lambda\in \Lambda}$ be an approximate unit of $A_\tmu$ and $\{w_\mu\}_{\mu\in M}$ one of $B_\tmu.$
    Consider $K=\Lambda\times M$ with the ordered $(\lambda,\mu)\leq (\lambda',\mu')$ iff $\lambda\leq \lambda'$ and $\mu\leq \mu'.$
    Given $\kappa=(\lambda,\mu)\in K$ set $d_\kappa:=\left(\begin{smallmatrix} v_\lambda & 0\\ 0 & w_\mu \end{smallmatrix} \right).$
    Then $\{d_\kappa\}_{\kappa\in K}$ is an approximate unit of $\Lb(\X)_t.$
    
    For $\kappa = (\lambda,\mu),$ $\xi = \left(\begin{smallmatrix} 0 & x\\ \widetilde{y} &0 \end{smallmatrix} \right)$ and $\eta = \left(\begin{smallmatrix} 0 & u\\ \widetilde{v} &0 \end{smallmatrix} \right)$ we have
    $$ \Lb(\gamma)_t(d_\kappa \xi )\eta =  \left(\begin{smallmatrix} \la \gamma_t(v_\lambda x),v\ra_l & 0\\ 0 &\la \gamma_t(yw_\mu),u\ra_r \end{smallmatrix} \right).$$
    
    By Cohen-Hewitt there exists $b,c\in B$ and $z,w\in \X$ such that $x=zb$ and $v=wc.$
    As $\be$ has a C*-globalization $\{\be_t(w_\mu b)c^*\}_{\mu\in M}$ converges to an element $p\in B_\tmu.$
    Then 
    \begin{align*}
	\la \gamma_t(v_\lambda x),v\ra_l
	    & = \la \gamma_t(v_\lambda a zb),wc\ra_l
		= \lim_\mu \lim_{\nu} \la \gamma_t(v_\lambda zw_\mu f_{\nu} b)c^*,w\ra_l \\
	    & = \lim_\mu \lim_\nu \la \gamma_t(v_\lambda zw_\mu) \be_t( f_{\nu} b)c^*,w\ra_l
		= \lim_\mu \la \gamma_t(v_\lambda zw_\mu) p,w\ra_l \\
	    & = \lim_\mu \la \gamma_t(v_\lambda zw_\mu \be_\tmu(p) ) ,w\ra_l
		= \la \gamma_t(v_\lambda z\be_\tmu(p)),w\ra_l.
    \end{align*}
    Note $z\be_\tmu(p)\in \X_\tmu,$ so that $\lim_\lambda \la \gamma_t(v_\lambda x),v\ra_l = \la \gamma_t(z\be_\tmu(p)),w\ra_l.$
    
    By symmetry $\{\la \gamma_t(yw_\mu),u\ra_r\}_{\mu\in M}$ is convergent. Then $\{ \Lb(\gamma)_t(d_\kappa \xi )\eta \}_{\kappa\in K}$ is convergent.
\end{proof}

\bibliographystyle{amsplain}
\bibliography{enveloping_2016_sent}

\end{document}